\documentclass[12pt,reqno]{amsart}
\usepackage{fullpage}
\usepackage{times}

\begin{document}

\title{Addressing the underrepresentation of women in mathematics conferences}
\author[Greg Martin]{Greg Martin} 
\address{Department of Mathematics \\ University of British Columbia \\ Room 121, 1984 Mathematics Road \\ Canada V6T 1Z2}
\email{gerg@math.ubc.ca}
\maketitle
\thispagestyle{empty}

\begin{abstract}
Despite significant improvements over the last few generations, the discipline of mathematics still counts a disproportionately small number of women among its practitioners. These women are underrepresented as conference speakers, even more so than the underrepresentation of women among PhD-earners as a whole. This underrepresentation is the result of implicit biases present within all of us, which cause us (on average) to perceive and treat women and men differently and unfairly. These mutually reinforcing biases begin in primary school, remain active through university study, and continue to oppose women's careers through their effects on hiring, evaluation, awarding of prizes, and inclusion in journal editorial boards and conference organization committees. Underrepresentation of women as conference speakers is a symptom of these biases, but it also serves to perpetuate them; therefore, addressing the inequity at conferences is valuable and necessary for countering this underrepresentation. We describe in detail the biases against women in mathematics, knowing that greater awareness of them leads to a better ability to mitigate them. Finally, we make explicit suggestions for organizing conferences in ways that are equitable for female mathematicians.
\end{abstract}
\bigskip

\section{Introduction}

In the context of mathematics conferences, the subject of gender is somewhat of a taboo. Certainly, bringing up the subject at all during a conference would be deemed outside the norm. And yet: in our graduate programs, women are still noticeably in the minority. We have a significant shortage of female mathematicians in our departments, particularly in senior positions. And where conferences are concerned, it is unfortunately quite common to have so few female speakers that they stand out from the homogeneously male pack. The pace of progress seems to be slowing, if we are continuing to improve at all. In short, we have too few women among the speakers at our conferences, and there is good reason to doubt whether the problem is simply going to fix itself any time soon.

The purpose of this article is to examine the issue of underrepresentation of women as speakers in mathematics conferences, as well as related disparities in conference organization, mathematics department composition, and recipients of mathematical prizes. We believe that it is our ethical responsibility to equitably represent all members of our profession\footnote{By ``our profession'', we are referring to academic mathematics; by ``our responsibility'', we mean the responsibility of all academic mathematicians, regardless of gender. This article is targeted most directly to American mathematicians, but it is relevant for mathematics around the world.} and to dismantle any obstacles to advancement in that profession, particularly when those obstacles disproportionately burden a minority group. In particular, we argue that it should be an explicit priority for any organizers of mathematics conferences to address appropriate representation of women in their lists of speakers; we further assert that we are not currently succeeding at meeting that priority.

%

Inviting speakers to conferences is about more than just rewarding a few already established people: we want to enrich the research of attendees and speakers alike. And one aspect of that enrichment is to expose ourselves to as many new and different viewpoints as possible; limiting our speaker pool (however unintentionally) is directly at odds with this goal. Research has shown that demographic diversity\footnote{In this article, we address gender diversity directly; however, much of what we say applies equal well to other aspects of diversity. We will include references to bias unrelated to gender, for example, when it illuminates the points we are making. Whenever we use the unmodified word ``diversity'', we intend the statement to be valid whether it is read as a statement about gender diversity specifically or about diversity in general.} has measurable positive effects on the outcomes of group enterprises; conversely, lack of diversity, in addition to perpetuating harmful stereotypes about mathematics, actually diminishes our ability to evaluate unfamiliar ideas. Once we include gender diversity among our explicit goals for conference organization, we become motivated to recognize the problem that currently exists and to proactively seek solutions.

Shortfalls of female speakers are, unfortunately, extremely common in all areas of science, technology, engineering, and mathematics (STEM).\footnote{Similarly, we address academic mathematics directly herein, but much of what we say applies equally well to other STEM fields and types of careers. We certainly make relevant arguments using data from STEM fields besides mathematics; indeed, we find corroborating evidence from disciplines even further afield.} In mathematics, just as in other STEM fields, American graduate schools have been producing a steady source of female PhDs for a generation, but that level of representation has failed to persist in many aspects of our discipline. We will display the overall underrepresentation of women, as well as additional differential underrepresentation in prestigious conferences roles, in two recent major conferences: the 2014 ICM and the 2014 Joint Meetings of the AMS/MAA. The same disappointing trends occur in employment statistics, on editorial boards, and in lists of prize-winners. Clearly, our current system is not living up to our standards where gender diversity is concerned; what are the causes of this shortfall?

Certainly there is no genetic predisposition that favors males over females in STEM fields (despite how often such claims are made). Girls and boys have always performed comparably on measures other than standardized tests; on these tests, the gap has dramatically decreased, to a nearly insignificant size, over the last generation. Furthermore, gaps in standardized test scores are significantly correlated to measures of gender inequality in the students' cultures. These effects manifest themselves in the set of high-achieving mathematics students as well as in the entire population (refuting ``more males at the top'' theories as well). None of this is consistent with innate gender-based differences in mathematics ability.

What, then, can be causing this underrepresentation of female mathematicians? It arises, in fact, from an assemblage of deeply entrenched biases that have been surreptitiously inserted into our perceptions and reactions. These implicit biases cause us to internally associate STEM careers---and, for that matter, positions of authority---with male defaults. Our culture's insertion of these biases into our subconscious, sadly, begins extremely early in our lives.

Through unconscious differences in the way they respond to girls and boys, schoolteachers reinforce assertive behavior and (over)confidence in boys, but passiveness and math anxiety in girls.\footnote{It is worth remarking that statements of this type are statistical statements, about large-scale trends of behavior. Of course there are individual exceptions to any such trend. However, unlike a proof in an axiomatic system, the existence of specific counterexamples does not invalidate the larger trends described in this article. Scientific research involving human behavior in complex societies does not look like theoretical mathematics, but it is completely appropriate for its subject.} Girls become less likely to volunteer for mathematical enrichment than boys, and children are unwittingly trained to perceive women as less mathematically able than men. These erroneous attitudes are compounded by a pervasive categorization of mathematical ability as being fixed and innate, rather than malleable and able to be strengthened---a categorization that our professional community unintentionally perpetuates.

Once these implicit biases are in place in all of us, they lead to further measurable discrimination that happens right under our noses. Experiments with dual versions of applications, resumes, and promotion files, identical except for the gender of the name, consistently demonstrate that women are rated lower than men even when there is literally no difference between them. Teaching evaluations display the same gender differential, as do multiple other evaluation instruments both inside and outside STEM; the vaguer and less concrete the evaluation criteria, the more our unconscious biases manifest.

The gender-related implicit biases of those around us also socialize us into choosing different behavior patterns. Without realizing it, we interpret a man's assertive demeanor as confidence but a woman's assertive demeanor as abrasiveness; we notice when women interrupt men but not when men interrupt women. The corresponding negative reinforcement indoctrinates women into dismissing their own ability (and socializes men into overestimating their own). Particularly in contexts that are stereotyped as male, such as STEM, this socialization of women leads to internal experiences of impostor phenomenon and stereotype threat, which degrade their ability to succeed. We cannot judge the choices that individual female scientists make in isolation from this social context.

The fact that biases exist at every stage of students' and professionals' careers causes a ``leaky pipeline'', where fewer and fewer women find success in advancing to progressively higher levels of achievement. Well-documented examples of this in the business world include the persistent wage gap between women and men and the poor record of top companies promoting women to the executive level. The same attrition can be seen in mathematics when looking at grant funding, tenure decisions, and awardee selection.

We are mostly unable to perceive, in individual situations, this pervasive pattern of invisible discrimination (that is the definition of invisible!); as a result, we fool ourselves into thinking that academic mathematics is a pure meritocracy. Once we take a closer look, however, we see that the current system (academic and societal) has been sullied with extraneous features that consistently discount merit where disadvantaged populations are concerned. Making an effort to address underrepresentation of women in mathematics, therefore, is not some extra component that introduces injustice---rather, it is an attempt to recognize and do away with the injustice that is already present.

It behooves us, then, to consider explicit actions we can take to mitigate the current unfairness in our discipline. Underrepresentation of women at conferences is a symptom of this unfairness, but it also contributes to perpetuating it; for this reason, we find it extremely important for this particular symptom to be treated (in conjunction with efforts to point out the larger inequities). Our goal is to put into practice guidelines for compensating for all the bias inherent in the system, with the hope that conscious attention to those biases will also help reduce them in the future.

To that end, we should plan our conferences with equitable gender representation in mind from the very start, and explicitly communicate with other conference organizers our expectation of meeting this goal. We should be extremely attentive to the way we select speakers, particularly keeping in mind that we are prone to misevaluating academic records of women and to overlooking qualified female candidates. We should recognize that mundane logistical choices can make conferences less welcoming to women if we are not careful. Moreover, we should publicly commit to equitable gender representation at our conferences, display this commitment visibly in conference materials and through our actions during the conference, and track over time how well we are (or are not) succeeding. Finally, we should simply talk more openly about underrepresentation of female mathematicians, not only in the context of conferences but in all aspects of the academic career; and we should have our words (and our attention to others' words) reflect the reality that mathematics is for women as well as for men.

We now go into more detail about the facts described in this introduction. We provide ample references to the primary research literature in sociology and psychology; we also include articles and informative summaries published by research organizations, web sites of relevant institutions, and well-written blog posts. In Section 2 we discuss the priorities for organizers of mathematics conferences to address, and we support the assertion that we are not currently succeeding with appropriate representation of women among speakers at conferences. We demonstrate in Section 3 that the underrepresentation of women in mathematics, far from being genetic in origin, in fact arises from an assemblage of deeply entrenched biases that have been surreptitiously inserted into our perceptions and reactions. With this perspective in mind, we conclude in Section 4 with lists of guidelines for conference organizers committed to equitable representation of female speakers.

\section{Recognizing our responsibility}

Dialogue about how to undertake any enterprise must start with a necessarily philosophical inspection of the priorities of that enterprise. We begin by discussing the various goals of organizing a conference, and in particular placing the priority of diversity itself---why it is important, what benefits it provides, and what harms are done by its absence---on equal footing with our other priorities. We go on to verify our empirical observation that this diversity priority is not being successfully met when it comes to women speaking at mathematics conferences. At the end of this section, we quickly review the research that refutes any attempt to ascribe these empirical observations to genetic differences between the sexes, in preparation for an examination of the actual causes of the disparity.

\subsection{Appropriate representation of women as an explicit priority}

Every project comes with a set of goals and priorities; organizing a conference and selecting its speakers is no different. Even though we tend to leave these goals unspoken, stating priorities always improves the quality and appropriateness of the end result. So let us articulate what we hope to accomplish when we choose a set of speakers for our conferences.

We know that giving a talk at a conference is not simply a reward for having published the most papers or won the most awards---otherwise every conference would feature the same prolific writers or the set of living Fields Medalists (or alternately, by induction, the speakers who have given the most talks at conferences!). Rather, the main purpose of a conference is to enrich the research and careers of those in attendance. In service to that goal, we try to include people working in newly hot topics, to give exposure to up-and-coming colleagues, and even to invite mathematicians who give accessible and entertaining talks~\cite{G}. Of course, inviting speakers is beneficial for the speakers' careers as well, and we value providing that benefit equitably to a representative cross-section of successful practitioners of our field. In addition, providing speaking opportunities to more mathematicians increases the community's awareness of them, leading to a greater likelihood of subsequent invitations to speak beyond the individual conference we are organizing. Thus the benefit extends those future conferences as well, by enlarging the pool of known candidates to invite as speakers or even to organize events around~\cite{AWM}.

In addition to these already laudable goals, much of what we want to accomplish in our conference organization involves diversity. We instinctively accept the desire for diversity when it comes to seeking a spectrum of mathematical fields, for example: a Joint Meeting of the AMS and MAA, or an International Congress of Mathematicians, would definitely appear strange if it had, say, no number theory whatsoever. Depending on the scope of our conference, we might also include as a goal the representation of scientists at various stages of their academic careers, and perhaps speakers from different geographical regions as well~\cite{VS}. In all these ways, we understand that diversity is good for its own sake, as the practical implementation of the desire to represent the entirety of our community. It stands to reason, then, that {\em appropriate representation of women as speakers is a valid and sensible priority for our conferences (indeed, for all STEM conferences).} Indeed, we deem it an even more crucial priority: while complementary thematic conferences can favor scientists from different topics and regions without any one population being worse off on average, gender differences consistently manifest as underrepresentation of the same population, namely women---a population, moreover, defined in terms of a characteristic completely unrelated to the subject matter of the conference.\footnote{Conferences specifically designed to be exclusively for female mathematicians exist, of course, but do not undermine this argument---indeed, such conferences are a specific reaction to the underrepresentation of women at mathematics conferences and the biases that female mathematicians face. If there were no gender bias in the discipline of mathematics, then organizing conferences for women only would not be necessary.}

It is worth pointing out that not only is diversity a noble ideal, but it also has measurable positive effects. Groups that are more diverse have been shown to benefit from the greater range of perspectives present and from a more inclusive environment that encourages people to contribute their individual views \cite{Su}. Including a diverse set of participants ensures that the best minds are represented, affirming the event's commitment to honoring merit---and merit-centered processes elevate the performance of everyone involved \cite{Ries}. With respect to gender, efforts to seek out submissions to tech conferences by women have resulted in noticeably more submissions at the highest end of the quality range \cite{JSConf}; as a related example, startup companies tend to be more successful when founded by women than when founded by men \cite{Ries}.

On the other side of the same coin, lack of diversity has measurable detrimental effects. Math, in our society, has plenty of negative stereotypes associated with it already; presenting nearly all-male conference speaker lineups perpetuates the stereotype of math as for men alone. This stereotype perpetuation occurs in other fields too: in philosophy, for instance, it ``undermines the self-confidence of women who aspire to become professional philosophers, or to remain in this exceptionally competitive profession [and] feeds the conscious or unconscious biases against women of the people who decide the fate of those who aspire to become or remain in the profession'' \cite{G}.

Psychological studies have shown that people unconsciously evaluate women less favorably in settings where they make up a small fraction of the participants, all the more so when gender-typing (the process in which our society trains us to associate certain activities or qualities with a single gender) is present \cite{H}. Thus an underrepresentation of women also causes them to be judged more harshly. The simple fact that a conference is overwhelmingly male makes it less welcoming for women to apply for and join (as it does for businesses \cite{Ries}), because the imbalance makes their gender a salient characteristic in a situation where it oughtn't be; the artificial salience of gender is a major component of ``stereotype threat'' \cite{FP:G}, about which we will say more later. Lack of gender diversity also harms men: homogeneous groups promote a sense of entitlement and complacency, an atmosphere where feedback and outside opinion is less welcome, and a susceptibility to ``groupthink''. As E.\ Ries summarizes \cite{Ries}: ``That's why I care a lot about diversity: not for its own sake, but because it is a source of strength for teams that have it, and a symptom of dysfunction for those that don't.''

With this list of priorities in hand, it is clear that we would not be satisfied with a conference whose speakers were all in stagnant mathematical fields, or uniformly dull and inaccessible, or (thematic conferences aside) all in their 20s or all from Texas. No organizer of a Joint Meetings or an ICM would say, ``Gosh, every time we plan our conference, we always end up with practically no number theory at all \dots\ oh well, them's the breaks!''---because an outcome that fails to achieve an explicit priority is not a matter to shrug off, but rather motivation to examine and improve the process leading to the outcome. Similarly, the reasonable reaction to chronic underrepresentation of women at math conferences---which is a failure to achieve an even more important explicit priority--- is not ``Gosh \dots\ them's the breaks'', but rather a concerted effort to identify and overcome the flaws in the process.

Of course, concentrating on gender carries with it the danger of failing to acknowledge other existing inequities present to various extents in the US and other countries (relating, for example, to ethnicity, sexual orientation, physical ability, socio-economic status, country of origin, religion, or political ideology); for that matter, we are starting to better understand that our society's binary gender construct does not fully represent every individual's self-identification (nor even their biology). However, we should not let the perfect be the enemy of the good.
There have been hundreds of math PhDs earned by women in North America alone, every year for many years running \cite[Table GS.2]{CMR}. We quite plainly are in a position to ensure appropriate female representation at math conferences; our inability to deal immediately with all inequities does not give us an excuse to ignore the important inequities we do have the power to correct. That being said, most of the discussion in this article can be immediately applied {\em mutatis mutandis} to related priorities such as better inclusion of underrepresented ethnicities.

\subsection{The present shortfall of female speakers} \label{miao}

Once we agree that proper representation by female scientists is a priority for us, it is depressingly easy to see that our priority is not being met. This failure is not unique to mathematics: even in scientific disciplines with greater gender parity, such as physical anthropology, primatology, and microbiology, women are less likely than men to be invited to speak, particularly when the organizing team does not include women \cite{IYH,CH}. In tech conferences, proposed talks are submitted far less often by women than by men, unless herculean efforts are made by the organizers to solicit proposals from women \cite{St}---even though the proposals that are eventually submitted by women tend to be of significantly higher quality on average than those submitted by men \cite{JSConf}.

Returning to mathematics, let us fix a conservative measuring standard and then apply it to two recent highly prestigious conferences. For this discussion, we will use 24\% as a minimum estimate of an appropriate proportion of female speakers in mathematical conferences. Every year since 1991, the percentage of PhDs in mathematics granted by US institutions to women has been 24\% or greater, with a peak as high as 34\%~\cite{AMSsurvey}. While these percentages would ideally be closer to 50\%, they still represent thousands of female mathematicians---mathematicians who have used these two decades and more to complete successful postdocs, amass tens of thousands of publications, earn tenure, be promoted to full professor, train PhD students of their own, and even lose their eligibility for the Fields Medal on account of their age. In particular, suggestions to the effect that not enough time has passed for women to work their way through the system cannot be taken seriously.

The first of the two conferences we examine is the 2014 International Congress of Mathematicians (ICM), held in August 2014 in Seoul, South Korea. Counting directly from the program on the ICM's official website, we find that only one of the twenty plenary speakers (5\%) was female. There were 20 sessions of invited speakers, 17 of which had at least six invited speakers; more than half of those 17 sessions included either a single woman among the speakers or no women at all. Overall, among the plenary and invited speakers, there were 35 female speakers out of a total of 237, or 14.8\%. (All the data gathered for this section can be examined in more detail in the author's supplement \cite{Martin} to this article.)

Human beings' notoriously poor sense of probability tempts us to believe that such underrepresentation of women might just be the result of chance. Indeed, it is a general trait known to psychologists that people tend to ascribe outcomes they expect to stable causes, while chalking up unexpected outcomes to temporary causes \cite{CI}. We mathematicians are well equipped, however, to perform the easy calculations showing how wrong this instinct would be. The appropriate null hypothesis is ``the ICM speakers were selected independently of gender from among the pool of people who have received PhDs in mathematics in the last 25 years''. Under our conservative 24\% assumption from above, the observation of nineteen male plenary speakers and one female plenary speaker rejects ($p<0.031$) this null hypothesis. Indeed, it is 18 times as likely that we would have seen an ``overrepresentation'' of female plenary speakers (five or more, since $20\times24\%=4.8$) by chance than to have seen at most one.

We might hope to discount this one observation as an anomaly (despite the fact that this statistical test is designed precisely to tell us not to do so); however, we have much more data at our disposal. The data on invited ICM speakers as a whole soundly rejects ($p<0.0004$) the null hypothesis that speakers were included independently of their genders; in this case, it is almost 1600 times as likely that chance would have led to more than 24\% women than at most the 14.8\% we observe. The only rational conclusion is that there has been bias somewhere along the line. One wonders, of course, where this bias takes place: in the process of producing the PhDs in the first place, or during postdoctoral positions, or faculty hiring and promotion, or evaluation of research records, or selection processes for conferences. One of the main purposes of this article is to demonstrate the existence of {\em gender biases at every one of these stages}.

The second conference we examine is the 2014 Joint Mathematics Meetings of the AMS and the MAA, held in January 2014 in Baltimore, Maryland. At first glance, the news seems to be better here: with 746 female participants (speakers and organizers) out of a total of 2,708, the percentage of female participants was 27.5\%, which at least exceeds our conservative measuring standard of 24\%. As we look a little more closely, however, we see some alarming disparities in the distribution of those participants. The percentage of speakers in contributed paper sessions who were female was 36.8\%, while the percentage of invited speakers who were female was only 25.6\%; this echoes the lower percentage of invited participants observed in other sciences \cite{IYH}. When we restrict to AMS-organized invited sessions, there were only 24.8\% female speakers and 22.0\% female organizers.

Among the JMM sessions that had their organizers listed explicitly, sessions with at least one female organizer had an average of 38.3\% female speakers, yet sessions with no female organizers had only 19.8\% female speakers, or half as many. Using a binomial regression, the null hypothesis (that the proportion of organizers who are female is not associated with the proportion of speakers who are female) is dramatically rejected: the probability of seeing a result at least as extreme as the data is less than $10^{-12}$. Again, this disparity echoes observations in other fields~\cite{CH}.

Differential representation of women can be found in employment statistics as well \cite{CMR}: as we pass from part-time faculty to full-time non-tenure-track faculty, to tenured faculty at non-PhD-granting institutions, to tenured faculty at PhD-granting institutions, the percentage of women decreases steadily. Where mathematics journals are concerned, a sample of ten of the most prestigious ones yields only 6.6\% of editors who are female; six of these ten journals have no women at all on their editorial boards. (Again, details can be found in \cite{Martin}.) And while it was legitimately exciting when M.\ Mirzakhani became the first woman to be awarded the Fields Medal in August 2014, it is hard not to wonder, given the fact that less than 2\% of all Fields Medalists and 0\% of all Abel Prize winners to date are female, how many outstanding female mathematicians have not had their work sufficiently recognized.

The AMS website has a page listing eighteen prizes and awards which have been given, as of the end of 2014, to 374 recipients; only 25 of them, or 6.7\%, have been women. But even that is not the complete story: one of these prizes is the Ruth Lyttle Satter Prize in Mathematics, for which only women are eligible. Once this prize is removed, the percentage of AMS prize and award recipients who are female plummets to a paltry 3.3\%. Rather embarrassingly for us, in the history of the AMS, more women (thirteen) have won the Satter Prize than the seventeen other AMS prizes and awards combined.

It is certainly not rare for us to experience firsthand the underrepresentation of women at conference after conference, award announcement after award announcement. The very fact that we can find, on the internet, a web page \cite{P} devoted to debunking the claim that underrepresentation of women at STEM conferences is often due to chance, as well as a ``bingo card'' \cite{F} full of anticipated excuses for not having more female speakers at an STEM conference, is dark-humor evidence that underrepresentation of women in science is both widespread and widely dismissed. Even single instances of underrepresentation we see are too skewed to rationally explain away as due to chance; in aggregate, there is simply no way to ignore the reality of bias. Despite the part of human nature that doesn't like to examine our imperfections closely, we simply must acknowledge the truth: the system we have in place now is causing a failure to meet our gender-diversity priority. And in addition to missing out on the positive effects of diversity, we need to question our default assumption that we really are choosing the most qualified speakers available. As Ries writes \cite{Ries}: ``[W]hen a team lacks diversity, that's a bad sign. What are the odds that the decisions that were made to create that team were really meritocratic?''

As scientists, after observing the existence of gender-based inequities in our discipline, our natural next step is to investigate their causes.

\subsection{Invalidity of genetic explanations}

A common interjection at this point in a discussion of women in STEM disciplines is to propose that there are neurological or ``hard-wired'' differences between female and male brains that result in men being better, on average, at mathematics than women. The simplest such hypothesis is that the distribution of males' mathematical ability (however that might be measured and quantified) is of similar shape to, but with a significantly higher mean than, the distribution of females' mathematical ability. A variant of this hypothesis is the just-so story that the two distributions have similar means, but the males' distribution has a significantly higher variance than the females' distribution; somehow the assertion that there are exceptionally math-hopeless men in the population seems soothing to those proposing that only men can be exceptionally math-savvy in great numbers.

Many researchers have made explicit various hypotheses that might explain overrepresentation of men in STEM fields \cite{Bar,EHL,KM,LHPL}, some genetic and some environmental. Any hypothesis, regardless of the motivations of its proponents, can stake its claim as a possibility for the truth and have its validity tested by the usual scientific method. As it happens, however, these genetic hypotheses have been thoroughly refuted by scientific studies. Indeed, such a refutation has appeared recently in the pages of these Notices \cite{KM}, so we restrict ourselves to a brief summary here.

First, the hypotheses that are genetic in nature are usually supported by girls' and boys' scores on standardized tests, even though restricting to this single instrument of measurement of children's mathematical skill is problematic \cite{AAUW,HLLEW,NV} and not particularly consistent with other data. (For example, girls consistently get better grades than boys in math classes at all levels of primary and secondary school \cite{EHL}.) It is true that a generation ago, boys did perform notably better than girls on standardized tests, but that gender gap has been repeatedly shown to be practically nonexistent today \cite{AGKM, EHL,GMSZ,HLLEW,HM,KM,OECD}; such a rapid change is utterly inconsistent with any genetic basis for differential mathematical skill. Moreover, several studies of students in dozens of countries \cite{EHL,GMSZ,HM} have shown that when gaps between boys' and girls' mathematical performance do exist, they are significantly correlated with measures of gender inequality in the country's employment, government, and culture.

The hypothesis that men exhibit greater variability in mathematical skill than women, and hence that men are appropriately represented in the subset of people with extremely high mathematical achievement, also does not stand up to scrutiny \cite{AGKM,GMSZ,HLLEW,HM,NV}. The correlations with measures of gender inequality mentioned above persist even when restricting to only the best mathematics students. It is worth noting that this ``greater male variability hypotheses'' would predict that the overrepresentation of men in scientific fields would be extremely large in the most math-intensive careers but would diminish as one examined jobs that required less and less mathematics for success; this prediction is certainly inconsistent with the empirical observation that underrepresentation of women is uniform over STEM disciplines and careers.

In response to a question to this effect at a panel discussion, N. d. Tyson \cite{Tyson} described the many ways in which his environment discouraged him, an African--American, from pursuing science as a career, and speculated that the experiences of girls interested in scientific opportunities might be quite similar. ``So before we start talking about genetic differences,'' he said, ``you gotta come up with a system where there's equal opportunity. {\em Then} we can have that conversation.'' Even though there seems little need, given the evidence summarized above, to continue the conversation about genetic differences, it is definitely true that there is no point in looking elsewhere if the system itself is not equitable to women and men. Accordingly, we turn now to an honest examination of the ways our society and our academic system treat women and men who are interested in mathematics.

\section{Assessing the obstacles}

We've debunked the arguments that there are genetic reasons for women to constitute such a small percentage of speakers at mathematics conferences. The opposing hypotheses assert, in one form or another, that the underrepresentation comes from environmental factors (aspects of our society, our academic system, and so on);  these hypotheses must clearly be taken seriously. As Ries writes \cite{Ries}, ``Demographic diversity is an indicator. It's a reasonable inference that a group that is homogeneous in appearance was probably chosen by a biased selector.'' 

To be clear, we are not suggesting that academic institutions in the US are prohibiting women from joining their universities or degree programs, as was the case several decades ago; we are not even suggesting that significant numbers of individual faculty members intentionally discriminate against women any more. ``The word `bias' here is not meant to imply deliberate bias. Although there may be deliberate cases, those are not the ones we are concerned about. Rather, we are concerned about the subtle, unintentional examples'' \cite{VS}. More specifically, we are proposing that typical mathematicians (indeed, typical human beings), both female and male, carry inside them unconscious habits and patterns of thought that add up to a significant difference in how female and male mathematicians are perceived and treated.

Of course, we are scientists: even though we have seen plentiful arguments against genetic hypotheses, we also need to see arguments supporting this new hypothesis. If we are to believe that the underrepresentation of women at our mathematics conferences is unknowingly caused by our community of well-meaning mathematics professionals as a whole, what are the mechanisms at work?

Fortunately for our argument, but most unfortunately for our discipline, the mechanisms causing bias against women in mathematics are both numerous and thoroughly documented.

\subsection{Ubiquity of implicit biases} \label{implicit bias sec}

We are complicated human beings in a complicated society, and as it turns out, we actually carry with us implicit biases hidden in the blind spots of our self-perception. Colored by these biases, our brains associate words like ``scientist'', ``boss'', ``doctor'', ``leader'', ``genius'' with male defaults. Mathematics, specifically, is stereotyped as a discipline for men, as has been measured by tests of people's implicit attitudes \cite{GRD,LHPL}. (We encourage the reader to try talking about colleagues or students without using gender pronouns: watch how quickly the assumptions of the listeners kick in.) These biases cause us to think of male colleagues more readily than female colleagues, ``leading to more invitations to men, leading to greater visibility for men, leading to yet easier availability of men's names'' \cite{VS}; they even cause us to evaluate files differently when the name is changed between female and male, with everything else kept equal. These biases, furthermore, are present in both women and men.

This sounds too terrible to be true (and indeed, those of us who have experienced success in our discipline are understandably resistant to the assertion of unfairness in the system). However, the biases summarized above are an extremely well-chronicled phenomenon.

\subsubsection*{Implicit biases in schools}

As early as primary school, our society trains girls to believe that they are not as valued as boys in an educational environment. Girls receive significantly less attention from teachers than boys \cite{AAUW}. Specifically, boys assert themselves by raising their hands, switching the arm they have raised, even jumping up and down until they get their teacher's attention, while girls will put their hand down soon when they are not called on. These behaviors are reinforced by the reactions they provoke: boys call out answers and contributions to the discussion out of turn and are not reprimanded, but when girls do the same, they are told to remember proper etiquette. When girls do answer questions, they are more likely to receive a brief response that merely recognizes the fact that they answered, whereas boys are more often given follow-up questions or time to expand upon their answer \cite{W}.\footnote{This unconscious differential treatment manifests with ethnicity as well as gender: for example, even though African--American girls attempt to participate more, they actually receive less attention than Caucasian girls from teachers \cite{AAUW}.}

With mathematics specifically, this situation is exacerbated by the fact that most primary school mathematics teachers are not math specialists; indeed, many primary school math teachers have math anxiety, as measured by psychological instruments such as surveys. Studies have shown that the math achievement scores of girls (but not boys) at the end of a school year are correlated with the level of math anxiety of their teacher, even though it was uncorrelated at the beginning of the school year \cite{BGRL}---presumably because witnessing their predominantly female teachers' anxiety reinforced the stereotype that boys are better at mathematics than girls. We are literally teaching girls to have math anxiety before they even leave elementary school.

Not surprisingly, these constant micro-inequities accumulate to form significant barriers. Girls assess their own mathematical abilities as worse than they actually are, in stark contrast to boys. Consequently, girls are less likely to volunteer for math competitions and programs for mathematically gifted children, and are more likely to underperform in standard testing environments, particularly when their gender is made salient \cite{NV}. While young children may not evaluate their peers' math ability differently according to gender, they do hold the attitude that adult men are better at math than adult women \cite{LHPL}. When rating their very own children, parents estimate the mathematical IQ of their sons as ten points higher on average than that of their daughters \cite{LHPL}.

A compounding factor is the dissonance between two theories of intelligence that people can implicitly hold: an ``entity theory'', which asserts that intelligence is ``a fixed quantity that cannot be changed very much by effort and learning'', and an ``incremental theory'', which asserts that intelligence ``is malleable and expandable'' \cite{MD}. In general, people who give credence to the entity theory (or who hear it asserted) are more likely to judge people's ability based on stereotypes than people who give credence to the incremental theory. Moreover, when comparing organizations or disciplines whose statements endorse fixed or malleable views of intelligence, people are more drawn to the latter kind and expect to feel more comfortable and accepted there \cite{MD}.

Believing math to be a fixed trait can ``turn students away from challenges that might undermine their belief that they have high ability'', while believing math to be a malleable quality can make students ``seek challenges that can result in better learning'' and ``remain highly strategic and effective in the face of setbacks, even showing enhanced motivation and performance'' \cite{GRD}. The entity theory, and the corresponding emphasis on innate ability, seems to be particularly strong in American culture: it has been shown to partially explain the difference in math achievement between students of American parents and students of Asian parents, who are more likely to subscribe to the incremental theory and emphasize effort and working through mistakes \cite{U}. The one-two stereotype punch of ``boys are better than girls at math'' and ``math ability is innate and fixed'' is actually what constantly erodes the confidence of female students as they develop mathematically; emphasizing a description of math ability as malleable and acquirable has been shown to mitigate the effect of the ``boys better than girls'' stereotype \cite{GRD}.

Not only do most of us mathematics instructors presumably hold an incremental theory of math ability (why else would we bother teaching at all?), we have surely all bemoaned many times the fact that our students can label themselves as innately incapable of understanding math, and wished we could make them understand that effort and persistence will pay off. But despite what we want our students to believe, our discipline seems to have a strong view of mathematical research ``power'' as a fixed aspect of intelligence. We gossip about our fields' wunderkinds and the {\it a priori} inevitability of their successes, we celebrate our Beautiful Minds and Good Will Huntings, and we prohibit mathematicians over 40 from winning the Fields medal on the suspicion that anyone capable of producing spectacular breakthroughs would have done so by then.\footnote{Whatever its origins or traditions, the rule prohibiting people over 40 years of age from receiving the Fields medal is another form of inequity---in this case, age discrimination. Again we note how the implicit (and questionable) assumption ``people over 40 are less capable of mathematical innovation than people under 40'' teams up with ``math ability is innate and fixed'' to produce this bias. In this case, the rule does not just exclude older mathematicians; it also causes a bias against people whose careers might have been delayed due to any number of reasons---parental leave, political unrest, bad job markets, taking care of elderly parents, settling into a tenured job later because of a two-body problem---unrelated to mathematical ability.} Furthermore, we fear making public mistakes (as evidenced by the comments of participants in the Polymath 8 project \cite{Polymath}) lest others judge us to be {\em inherently} incapable of serious mathematics.

In reality, our mathematical heroines and heroes, past and present, have all spent huge amounts of time learning new things, struggling to understand them, and feeling proud when they sort the dissonances out, just like us. It is important to explicitly remind ourselves of the acquirable nature of math research ability---particularly since we are not just instructors, but also personnel evaluators, and thus gatekeepers of the academic mathematical world.

\subsubsection*{Implicit biases in academia and business}

After receiving an entire childhood of unintentional training in the supposed differences between female and male math ability (and other cultural myths), we enter our careers burdened with implicit biases that manifest in a whole new set of ways. For example, professors at US universities who are contacted by students interested in their doctoral program respond more frequently to men than to women (and, for that matter, more frequently to Caucasians than to applicants of other ethnicities)---and this propensity is exaggerated in more lucrative fields and at more prestigious institutions \cite{MAC}. Both female and male faculty members rate students' application materials differently when the applicant is female or male: even with identical files, the female applicant is judged to be less competent, and male applicants are offered a 14\% higher starting salary and more mentoring on average than female applicants \cite{MDBGH}.

Similar disparities are present in situations other than authority figures evaluating younger personnel. Teaching evaluations of female professors by students are lower than teaching evaluations of male professors---all the more so when the professor teaches a technical course or requires a large amount of work from the students, or when the student evaluator has been given negative feedback by the professor or is male \cite{Kas,N,SK}. Topics for conferences are often chosen based on the strengths and specialties of the first star researchers who come to the organizers' mind---which are overwhelmingly male---and so even the very choice of conference topics can suffer from a bias towards men over women \cite{VS}.

In addition to these and many other examples in STEM fields (see also \cite[Chapter 8]{DK} and \cite{Newsome}), biases of this sort manifest in nearly every arena imaginable. In law school, traditional pedagogical structure depresses the grades of female students, causing a gender gap that disappears when the structure of courses changes~\cite{HK}. Orchestras started hiring dramatically more women when they began to place auditionees behind screens so that their gender could not be ascertained \cite{GR}. Prospective investors listening to entrepreneurial pitches pledge more funding to male presenters than female presenters, even when the content is the same \cite{BHKM}. Particularly in stereotypically male arenas, successful women are liked less and belittled more than successful men, to the detriment of their evaluations and careers \cite{HWFT}. Female parents are perceived as being less competent in the workplace and are given lower salaries than male parents \cite{DK}; managers display bias against women's requests for flexible work schedules, interpreting them for example as revealing less dedication to their career, unlike similar requests from men \cite{BGS,Munsch,RP}. 

A.\ Phillips (as quoted in \cite{G}) speaks of ``a cluster of vaguer characteristics which can override the stricter numerical hierarchy of grades or publications or degrees, always moderated by additional criteria. These more qualitative criteria [(]`personality', `character', whether the candidates will `fit in') often favour those who are most like the people conducting the interview: more starkly, they often favour the men.'' In other words, our biases find ways to manifest even when we are addressing objective, quantifiable data related to people. For example, observers tasked with judging people's attributes such as height, weight, and (startlingly) income from photographs consistently overestimate these attributes in men and underestimate them in women, even when objective comparison tools (such as a door frame of fixed height) is present \cite{BMN}. More saliently, when evaluating the research records of female and male scientists by their number of publications and the journals in which they appear, evaluators devalued the work of the women to the extent that a woman's file had to contain $2.5$ times as much productivity and impact as a man's file for the woman to be considered as competent as the man \cite{WW}. Sadly, the more subjective our criteria for research excellence are, the more our unconscious biases manifest and skew our evaluations.

Although we are most concerned with bringing unconscious sexism into the light where it can be examined and addressed, we must point out that explicit sexism is sadly not absent from STEM communities. Sexual harassment is so common at tech conferences, unfortunately, as to make it necessary to include explicit anti-harassment statements in conference materials \cite{GF:C}; such harassment is present in academic departments as well \cite[Fostering Success for Women in Science and Engineering, pages 6--8]{WISELI:o} and, in horrifyingly graphic form, on the internet~\cite{Hess}. B.\ Barres, a transgender mathematician who was born genetically female and began identifying as male to his colleagues during graduate school, reports \cite{Bar} that when he was outwardly female, professors would not give him credit for solving difficult math problems, fellowships were given to less-qualified male applicants, and he was frequently interrupted in conversation; one person even declared his research after he became outwardly male to be ``much better than his sister's''.

In summary: we don't want to have implicit biases, but, rather demoralizingly, we all do. This realization surely motivates us to consider what steps we can take to counteract our biases (and Section~\ref{process} contains many explicit suggestions). But before we get there, we need to further examine the outcomes that these implicit biases generate in our society and in our profession.

\subsection{Gender-based socialization differences}

We have seen how our perceptions of other people are tainted by implicit biases; but even our perceptions of ourselves are not immune. Much of our own behavior and self-evaluation is influenced by socialization---the inductive training we receive, from how members of our society typically react to certain actions from certain types of people, in how to act in the way society deems acceptable. Girls, for example, are socialized from an early age to place others' needs over their own interests, and men are socialized to expect women to act that way. Women who violate this social norm are deemed demanding and malicious instead of ``nice'' and likable \cite{BLGS,BBL,EK,HWFT}; ``what appears assertive, self-confident, or entrepreneurial in a man often looks abrasive, arrogant, or self-promoting in a woman'' \cite{EIK}, and indeed the word ``abrasive'' itself seems to be {\em de facto} reserved for women in performance reviews \cite{Snyder}.

Analogously to how boys and girls are (as described earlier) socialized in primary school classrooms, adult men speak more often and more forcefully, adopt more dominant body language, and interrupt other speakers more often than women. These habits of assertiveness are negatively reinforced in women, who are again disliked and perceived as untrustworthy when they adopt such patterns \cite{R}. When women and men give feedback to others, ``the evaluation of women depended more on the favorability of the feedback they provided than was the case for men'', and women (but not men) who gave negative feedback were judged less competent by the people they criticized~\cite{SK}.

Differential socialization of women and men means that corresponding behaviors are selectively reinforced. When conferences require submission of proposals, for instance, women are quick to dismiss their own ability to submit high-quality proposals, tending to submit in much smaller numbers in the absence of specific invitations. (We scientists detect patterns, after all: when a conference has startlingly few women year after year, and when women face greater scrutiny and criticism than men for equivalent work, it's logical for a woman to conclude that it might not be worth applying to speak!) On the other hand, men are very quick to submit---even when the quality of their proposals is well below average for the conference---because they overestimate their own abilities \cite{JSConf,St}. Analogously, in the business world, even an even-handed boss will end up making decisions that are biased against women due to not taking men's propensity for bragging about their achievements into account \cite{RSZ}.\footnote{This sentence, containing the word ``boss'', probably evoked a quick image of a person in the reader's mind. What gender was that person?}

Many women have ``internalized into a self-stereotype the societal sex-role stereotype that they are not considered competent'' \cite{CI}. For instance, ``[g]irls and boys with the same math test scores have very different assessments of their relative ability\dots. Conditional on math performance, boys are more overconfident than girls, and this gender gap is greatest among gifted children''; and socialized differences, such as men's overconfidence in their own abilities and women's reluctance to enter into competitive situations, are made more extreme when women are inadequately represented in the situation in question \cite{NV}. ``In line with their lower expectancies, women tend to attribute their successes to temporary causes, such as luck or effort, in contrast to men who are much more likely to attribute their successes to the internal, stable factor of ability'' \cite{CI}.

In fact, the aggregate effects of these socialization biases are so powerful that they manifest in a psychologically significant way. The phrase ``impostor phenomenon'' was used by Clance and Imes \cite{CI} to ``designate an internal experience of intellectual [phoniness], which appears to be particularly prevalent and intense among a select sample of high achieving women. ... Despite outstanding academic and professional accomplishments, women who experience the imposter phenomenon persists [sic] in believing that they are really not bright and have fooled anyone who thinks otherwise.'' This phenomenon affects a large and diverse group of women with fantastic career accomplishments \cite{Kap}.

Another internalized obstacle to success for women (and other minorities) is ``stereotype threat'', described by Spencer, Steele, and Quinn \cite{SSQ} in this way: ``When women perform math, unlike men, they risk being judged by the negative stereotype that women have weaker math ability. We call this predicament stereotype threat and hypothesize that the apprehension it causes may disrupt women's math performance.'' The effect of stereotype threat on actual measurable performance has been pointed out multiple times: for example, ``ability-impugning stereotypes such as these can trigger psychological processes that can undermine the performance of stereotyped individuals, including females in math'' \cite{GRD}. It has even been shown that when administering math tests, the likelihood of a girl ``choking'' on the test is noticeably diminished when the test is explicitly described in advance as not having displayed any gender difference in performance \cite{NV}; this effect has been observed in girls as young as five years old \cite{LHPL}. Seemingly innocuous environmental cues, such as references to words like ``pink'' or ``Barbie'', can prime stereotypic beliefs to the point where it changes what actions people deem acceptable, thus breeding more biased beliefs about themselves and more biased responses from others \cite{UC}.

We all have our behaviors affected by these socializations, but we certainly didn't choose them: they were inculcated into us by the culture in which we were raised, without our consent. If our society irrationally taught everyone continuously that blue-eyed people are less intelligent than other people, it would be the blue-eyed who experienced the impostor phenomenon and suffered from stereotype threat. And it would be foolish to think that this society's problems could be fixed just by blue-eyed people deciding to act or think differently: destroying the ``blue-eyed bias'' would be the responsibility of everybody whose default attitudes had been trained by this society---that is, each individual person, regardless of eye color. Similarly, in our current American society, with its (non-hypothetical) gender-based implicit biases, it is naive to think that women can simply change the way they react to their environment and dissolve all this inequity themselves. There are deeply entrenched reasons why resolving our underrepresentation problem will not be possible until all of us, not just the affected population, decide to devote effort towards recognizing and addressing the causes.

\subsection{Vicious circle of underrepresentation}

Without such efforts, the cumulative effect of biases throughout the system leads to the observation known as the ``leaky pipeline'': the higher the academic rank, the smaller the percentage of women (see \cite{CMR} and \cite[Advancing women in science and engineering]{WISELI:o}). We have already discussed how this leaky pipeline cannot be attributed to a shortfall of time for women to gain seniority in the academic world; the result must be due to aspects of our system that are biased against women. Inequities of this type have been categorized in sociology literature on the ``theory of cumulative advantage''~\cite{DE}.

In business, the fact that successful women are disliked and viewed as less competent than equally successful men prevents women from advancing---in rank, salary, and authority---as quickly as men \cite{HWFT}. For example, women are less likely to negotiate than men for raises and promotions; but when they do, the culture of many companies punishes them \cite{BLGS,BBL}, leading to a cycle of wage discrimination which, although having decreased somewhat over time, still persists today in the form of a significant gender wage gap \cite{EK,HWHH}. In addition to being implicitly regarded as unsuitable for leadership positions, women have fewer female contacts in positions of authority, which means that they are disadvantaged by having less influential networks \cite{EIK}. When they do actually become managers and leaders, ``[t]he mismatch between qualities attributed to women and qualities thought necessary for leadership places women leaders in a double bind and subjects them to a double standard'' \cite{EIK}; they have higher expectations placed upon them yet are perceived as less competent, and thus ``a woman manager's efforts to assert authority over others is subtly undercut by continuing, implicit assumptions that she is not quite as competent in the role as a man would be'' \cite{R}. Performance evaluations for women are far more likely to contain critical comments, and overwhelmingly more likely to have negative personality criticism attached to the comments on their performance, than those for men \cite{Snyder}. As a result, the leaky pipeline manifests itself in, say, disappointingly little change in the proportion of CEOs of top companies who are female \cite{CE}.

We have already seen biases when evaluating younger members of our discipline \cite{MAC,MDBGH,RSZ,WW} and when evaluating instructors \cite{Kas,N,SK}, both of which can obviously suppress the affected population from advancing in their careers; unsurprisingly, these biases manifest even when peers are appraising peers. Experiments in which identical conference abstracts were attributed to women or men show that peers (of both genders) perceived higher scientific merit, and were more likely to want to collaborate, when a male name was attached. Similarly, actual recommendation letters of contemporaries in STEM fields have been shown to exhibit different patterns of language usage in a way that benefits men over women \cite{KGH}. In medicine, for instance, letters of recommendation written for female faculty, in comparison to those written for male faculty, are shorter, mention high-status terms less and doubt-raising phrases more, lack basic features of recommendation letters more frequently, and tend to depict women as students and teachers while depicting men as professionals and researchers \cite{TP}. Even when a woman's tenure file is evaluated positively, the evaluators are four times as likely to volunteer ``cautionary comments'', saying that they would need to be given additional information to make a final judgment, than for a man's file \cite{SAR}. Once the leaky pipeline gets going, the mere fact of lower representation of women can exacerbate these biases: when MBA students assessed applicants for managerial positions, ``personnel decisions [by] both males and females were significantly more unfavorable towards women when they represented 25\% or less of the total pool''  \cite{H}.

These biased aspects of evaluation further disadvantage women in award competitions, leading to future disparities that, naturally, measurably compound themselves thereafter. For example, biased evaluations lead to smaller grants for women, which lead to somewhat curtailed research opportunities, which lead to artificially diminished research records that penalize them further for the next grant applications \cite{CGFSH,KGH}. Tenure cases are more harshly judged for female professors than for male professors \cite{SAR}. Male faculty members, who have benefited from these systemic biases enough to be overrepresented in elite universities, tend to take on fewer female graduate students and postdocs; this bias leads to fewer women being trained at these elite universities, which leads to female applicants being underrepresented in the next round of faculty hiring \cite{SS}.

The fact of the matter is: when a woman achieves the same level of accomplishment as a man (even if we were skilled at deciding who really was at ``the same level''), it is most probably a sign of much higher skill and potential in the woman, due to the cumulative disadvantages impairing her progress at every stage.

\subsection{The status quo as a distortion}

The above examples, and plenty of others from the literature, put the lie to the pleasant fiction that academia is a pure meritocracy where rewards always correspond to ability and achievement and nothing else \cite{G,KGH}. Instead, implicit biases and gender-based socialization sustain a persistent pattern of invisible discrimination against female scientists. That being said, each one of us might hope that we, personally, are objective and free from biases of this sort; unfortunately, that is extremely unlikely to be the case. Studies definitely document a universal tendency for individuals to believe themselves much less biased than others; and when our biases are hidden in our blind spots, our human nature comes up with all kinds of rationalizations about why the unsuccessful outcomes were not under our control. As L.\ Bacon writes \cite{Bac}: ``The defining feature of a blind spot is that we don't know it's there. And it's hard to notice it until we're challenged on it. We see this again and again with all-male speaker lineups at tech conferences. I certainly don't believe the organizers of those conferences are rabid misogynists; they just have a blind spot when it comes to gender, and frequently don't notice the lack of women until it's pointed out to them.''

The fact that we are aware of our thoughts and introspections yet do not recognize subliminal biases among them, psychologists theorize, convinces us that they must be absent even when they are surely present \cite{PGR}. And it takes courage, certainly, to look inwards and acknowledge our own imperfection (in this or any arena). It can be an eye-opening experience to take an implicit biases test \cite{GBN} and see that we are not a perfect specimen of objectivity. Indeed, being aware of our personal biases is far superior to the alternative, since people who consider themselves extremely objective can actually be more prone to act in a biased fashion \cite{UC}.

Even more unfairly, our biases against women operate on multiple levels in society as a whole. Even before getting to the point where science can display its biases and business can discriminate with its wages, girls have to deal with social pressures pushing them away from cerebral subjects, athletic pursuits, and male-typed careers; they grow up learning a version of history and science consisting almost entirely of the contributions of white males; they are expected to contribute more towards childcare, home maintenance, and so on outside of work yet are not credited for this ``second shift'' \cite{HoMa}; they also, unfortunately, can never ignore the very real physical danger of sexual assault \cite{AAUW,AGKM,BMN,LHPL}. ``What is interesting about the age old gender system in Western society,'' write C.\ Ridgeway and S.\ Correll \cite{RC}, ``is not that it never changes but that it sustains itself by continually redefining who men and women are and what they do while preserving the fundamental assumption that whatever the differences are, on balance, they imply that men are rightly more powerful.'' And it is difficult and frustrating to engage with this demoralizing reality. Not engaging with it, however, is tantamount to perpetuating it. In an ideal world, every part of our society would be gender-blind when it came to opportunities and rewards; but in today's actual society, being gender-blind effectively means ``if the status quo is biased, then we're ignoring that''.

Regarding mathematics conferences in particular, organizers might bemoan that there was ``nothing more'' they could have done to increase the number of female speakers, which suggests that our current system is neutral and that efforts to include adequate representation of women would be add-ons of some sort. But in fact, we have documented how our current system already includes a powerful collection of add-ons that consistently skew our decisions unfairly in favor of male mathematicians over female mathematicians. In other words, we are not simply trying to react to a perceived shortage of female mathematicians---we are, unintentionally and against our wishes, maintaining the shortage ourselves. So let us frame the issue of appropriate representation of women in mathematics, not in terms of some additional constraints that we must add in, but rather in terms of how to take out (or at least circumvent) the extraneous biases that are already there.

\section{Improving our process} \label{process}

As pointed out earlier, we would not accept consistently failing to meet other priorities, such as the representation of all mathematical areas, or inclusion of speakers from all geographical regions represented by the event: if a major conference's scientific structure were causing such a failure, we would augment or adapt its structure in response. While we might prefer to ``make it up as we go'' each time, the simple fact (as we saw near the end of Section~\ref{implicit bias sec}) is that such a strategy is inherently biased against success. Once we assert gender diversity as one of our priorities and acknowledge that the current system hinders our ability to fulfill that priority, it becomes clear that we must consider an explicit process for our conferences to assure adequate representation of female speakers.

We have gathered in this section some guidelines for meeting this priority. Many of them are common sense, especially now that we understand the causes and ubiquity of underrepresentation. While some of these steps might be labeled drastic if adequate representation of women at our conferences were in the ``Gosh \dots\ them's the breaks'' category, they are all perfectly natural ways of addressing a recognized failure to meet an important goal.\footnote{In a hypothetical situation where a national conference repeatedly left out a subject area such as number theory, the analogues of these steps would be natural ways of ensuring that number theorists were included. We reiterate, however, that remedying gender inequity is more important than such a hypothetical situation: while number theorists might be occasionally undervalued, female mathematicians are {\em consistently} underrepresented, based on a factor completely unrelated to the practice of mathematics.} While many of these action items will be useful no matter how the conference is structured, some of them make more sense for invited speakers, others for contributed lectures.

We repeat that our suggestions are aimed at mitigating existing unfairness in the current system of conference organization (which none of us desire to be there). This unfairness can be reduced and eventually eliminated, both by taking deliberate steps to fully include women in our scientific activities and by focusing attention critically on the unfairness itself. Moreover, in addition to resulting in appropriate inclusion of female mathematicians, we believe that thoughtful adoption of these guidelines will quite simply lead to better conferences, independently of speakers' genders.

\subsection{Plan from the beginning}

\begin{itemize}
\item The conference's scientific committee must have women adequately represented. This is a no-brainer. A group that cannot find women at this early, controlled stage of the process will not magically be able to balance genders at their conference later on. Recall our observations about the correlation between female organizers and female speakers \cite{CH,Martin}.
\item Plan far ahead of time, so that women, who tend to have more non-work responsibilities \cite{HoMa}, have adequate time to make arrangements for travel.
\item Don't automatically structure a conference (or part of it) around an eminent man, but consider building it around a woman---her expertise, her dates of availability.
\item Some conferences are built around a slate of invited speakers, while others consist predominantly of speakers selected from submitted proposals. The challenges to achieving gender diversity in a conferences are somewhat different for these two types of conferences. Be aware of, and plan to address, the challenges particular to the chosen format.
\end{itemize}

\subsection{Communicate expectations and goals among planners}

\begin{itemize}
\item As a planning committee, articulate all the priorities of the conference, including proper representation of women.
\item Set explicit targets---for example, that 30\% of speakers should be female. Note the subtle but important difference between ``targets'', which connotes the quantification of a goal we already wish to achieve, versus ``quotas'', which carries overtones of the myth of excluding more qualified men. Remember that having fewer women participate means that those who are present are judged more harshly \cite{H}.
\item Make the selection criteria explicit and unambiguous within the scientific committee, to minimize the effect of biases that give more credence to candidates who conform to irrelevant stereotypes.
\item Explicitly inform all decision-makers that letters of recommendation, and statements of self-promotion, tend to bias the process against women. A detailed list of suggestions and practices can be found in \cite{Valian}.
\end{itemize}

\subsection{Come up with names}

\begin{itemize}
\item Be explicitly prepared for the fact that men's names will come to mind more readily then women's names.
\item Do a literature search for women's names in the relevant areas. For that matter, make sure the relevant areas are broad enough to capture a significant cross-section of scientists.
\item Consult a database of female mathematicians. For example, if a conference is affiliated with any of the mathematical institutes sponsored in part by the National Science Foundation, the organizers should have access to their expertise database \cite{NSF}.
\item For that matter, organizers of recurring specialist conferences can band together with others in the field to build their own databases of minority scientists, such as the Women in Number Theory database~\cite{WIN}.
\item Have each organizer ask five respected colleagues to suggest female speakers for the conference, as part of the initial planning process.
\item Since women are overrepresented at lower prestige institutions, don't stop searching once the people at high-prestige institutions have been exhausted. Don't dismiss possible presenters because of quick initial impressions that they aren't on a high enough level---remember that implicit biases color these impressions.
\item Proactively seek personal communication with potential female speakers, brainstorming their possible submissions.
\item Add a ``suggest a speaker'' form or email link to the conference proposal website. Include nearby a link to a conference diversity statement (see below).
\end{itemize}

\subsection{Select speakers attentively}

\begin{itemize}
\item Give scientific panels and submission evaluators ample time to consider their decisions. Include live conversations in those deliberations, in which diversity priorities are discussed.
\item Structure the submission process so that proposal materials will be evaluated without seeing the submitter's gender.
\item Invite the entire first round of speakers, with women proportionally represented, at once---don't invite a round of speakers and then add a few last-minute invitations to women later.
\item Check that diversity is coming in prestigious conference positions (plenary talks, invited talks) rather than only contributed talks and poster sessions.
\item Communicate individually with invitees, rather than through form letters.
\end{itemize}

\subsection{Create equitable logistics}

\begin{itemize}
\item Set aside some of the budget specifically to offer to defray the costs of female speakers---after all, they have already been implicitly biased against for funding opportunities so far in their careers. Indeed, apply for grants from scientific funding agencies specifically for these costs, or partner with mathematical institutes or hosting institutions to acquire such funding. The Association for Women in Mathematics (AWM) advertises its own funding source, the Travel Grants for Women Researchers, as well \cite{AWM}.
\item For those in a position to grant funding to organized scientific events, make it a condition of acceptance that the organizers have an explicit process for ensuring proper representation of women at their event.
\item Be aware that travel logistics can be extremely difficult and stressful for parents of young children. Make arrangements for childcare and nursing mothers at conferences, and communicate that information to prospective participants. The official statement of the Association for Women in Mathematics on childcare \cite{Sa} is a good place to start. Follow the lead of the 2015 Joint Meetings of the AMS and MAA \cite{JMM}, for instance, or the case studies of conferences offering childcare arrangements in \cite[Chapters 5--6]{DK} and \cite{FP:C}. Particularly if one has never nursed an infant before, one should consider firsthand suggestions for how to make conference welcoming for nursing mothers~\cite{Lalin}.
\item Pay attention to ``extracurricular'' details, so that all aspects of the conference are both inviting and safe to women. For example, make sure that women's restrooms are as convenient and well-maintained as men's washrooms. As Geek Feminism points out \cite{GF:T}, ``conference dinners with 90\% or more men and free alcohol are not welcoming or safe''---a fact that many men, with the privilege of not having to worry about harassment, might overlook. Similarly, going out for evening drinks, at a location far from the conference venue or hotel, is simply not comfortable for a woman when the rest of the group consists only of men. Having to avoid social opportunities because of unsafe circumstances means missing mathematical discussions and networking opportunities, which has real professional consequences.
\end{itemize}

\subsection{Walk the walk}

\begin{itemize}
\item Track the diversity of a conference's speakers. Have the results publicly accessible.
\item If there has been lack of diversity, step forward and admit it---and make a public pledge to do better.
\item Have a diversity statement (such as \cite{OM}) prominently featured in conference materials, including the web site, and include it in calls for speakers and communications with potential participants. Even aside from helping to generate more female speakers, such statements can help women feel more comfortable attending the conference \cite{Ries,RM}.
\item As a conference organizer, attend sessions by female speakers (as well as speakers from underrepresented minorities), and initiate conversations with female participants about their research.
\item As a chair of a conference session, call on women for questions that follow talks as well as men (and notice how men speak without being called on more often than women).
\item Also include an anti-harassment statement or code of conduct (such as \cite{GF:C}) in conference materials.
\item Particularly for those who are part of a population with the privilege of rarely being the target of discrimination or harassment, proactively seek information about past experiences by surveying conference attendees, with such questions as: when you were deciding whether to attend this conference, what factors weighed in your decision? How inclusive do you find this conference in terms of gender, ethnicity, sexual orientation, people with disabilities, etc.\ for both speakers and participants? Have you experienced any harassment at this conference, including harassment based on gender, sexual orientation, or ethnicity?
\end{itemize}

\subsection{Talk the talk}

\begin{itemize}
\item Talk openly about underrepresentation of women at mathematics conferences (and in mathematics in general). Make discussing the issue, and whether our community is adequately addressing it, the norm, so that ignoring the issue becomes the controversial stance instead of the safe option.
\item Introduce good management practices, such as those described in \cite{Keyfitz}, for equitable hiring, retention, and work environment in the mathematics department.
\item Put on professional web pages a statement (such as \cite{VS}) endorsing diversity, or pledging participation only in events that make diversity an explicit priority.
\item Thoughtfully monitor interactions with other mathematicians. Do we mention the physical appearance of some graduate students but not others? Do we ask some postdocs about their research but others merely who their supervisor was? Do we interrupt some colleagues in conversation but not others? Are there any correlations between these differences and the other mathematicians' genders?
\item Practice representing mathematics as suitable for girls and boys, for women and men. Practice using gender-neutral speech patterns when speaking about mathematicians.
\item Speak about mathematics skill as a malleable quality rather than a fixed quantity \cite{GRD,MD,U}---not just to students, but among ourselves.
\end{itemize}

\section{Considering the issue further}

We all know that excellence in mathematics requires hard work, mental focus, and self-awareness, as well as an understanding of what details to focus those virtues upon. Excellence in conference organization is no different. We accept the necessity of this effort, both in research and in interactions with our scientific community, because we understand both the value of the end product and how introspection and attention to detail improves that end product.

It is an unfortunate reality that mathematics still has a gender inequity problem, despite the improvements we have made over the generations. There is good news, though: not only do we understand quite a bit about what causes underrepresentation of woman---as well as actions we can take to rectify it---but also, we can make our discipline's challenges easier to successfully address simply by talking openly about them.

For those wanting to learn more about the facts and guidelines included in this article, its list of references has been expanded into an annotated bibliography \cite{Martin} with remarks and quotes to help describe the content of those sources. We emphasize here several sources rich in further references to the relevant literature on the disappearing gender gap on standardized tests \cite{AGKM,EHL,GMSZ,KM}, biases and barriers to advancement for women \cite{DE,EIK,KGH,N,R,WISELI:o}, psychology theory and research \cite{EHL,GRD,PGR}, and action items to address gender bias in schools \cite{AAUW,W} and in academia \cite{DK,FP:G,WISELI:o}.

We have made the conscious choice to include only initials and last names in the bibliography and throughout this article. We have observed a tendency to be curious about the gender of the authors of the research we cite, and perhaps to involuntarily wonder how the authors' genders should affect our evaluation of their conclusions. These reflexive speculations, we believe, tellingly illuminate the depth to which these implicit biases about gender are ingrained in us, even though we rationally know that possessing one gender or another does not affect a person's objectivity.

We hope this article has been thought-provoking in ways that we as a community will carry with us into future planning. Only by examining our current practices, and carrying the resulting knowledge forward, will our conferences (and our departments) improve in this regard. Being socialized to have biases is not our fault; but preventing our biases from negatively affecting the world around us is nonetheless our responsibility.

We conclude with a challenge to each reader: find a female colleague who is willing to donate some of her time, and ask her about her experiences as a scientist in training and as a woman in today's society. Female readers will probably find commonalities of experience, while male readers might well be surprised at the injustice some mathematicians have had to deal with. The more concretely we apprehend the inequities that still exist, the better equipped we will be to remove them.

And finally, let us completely reverse the taboo: let us make it the norm to talk about gender in conferences, so that overlooking the issue becomes a dubious exception. None of us want to perpetuate prejudice, but we are doing so nonetheless. Let's change that.

\section*{Acknowledgments}

We thank J. Bryan for statistical analysis of the data from the 2014 JMM and W. Miao for gathering the data in Section~\ref{miao} and for locating several of the papers cited herein. We also thank T. Gowers and L. Walling for helping to crystalize our ideas at the early stages of conceiving this article, and L. Addario--Berry, J. Gordon, T. Gowers, C. Hagan, E. Jones, M. Lalin, W. Miao, A. Mottahed, R. Pries, and R. Scheidler for their helpful feedback on a draft of the manuscript. We are also grateful to other friends and colleagues, too numerous to list here, for their encouragement and inspiration.


\begin{thebibliography}{99}

\bibitem{AAUW}
American Association of University Women, How schools shortchange girls: executive summary, The AAUW Report, 1992. http://www.aauw.org/files/2013/02/how-schools-shortchange-girls-executive-summary.pdf (accessed September 15, 2014)

\bibitem{AMSsurvey}
American Mathematical Association, Annual survey of the mathematical sciences: full reports. http://www.ams.org/profession/data/annual-survey/survey-reports (accessed September 5, 2014)

\bibitem{AWM}
Association for Women in Mathematics, Mathematics travel grants. \\
http://sites.google.com/site/awmmath/programs/travel-grants/mathematics-travel-grants (accessed August 15, 2014)

\bibitem{AGKM}
T. Andreescu, J. Gallian, J. M. Kane, and J. E. Mertz, Cross-cultural analysis of students with exceptional talent in mathematical problem solving, Notices of the AMS 55 (2008), 1248--1260.

\bibitem{BLGS}1
L. Babcock, S. Laschever, M. Gelfand, and D. Small, Nice girls don't ask, Harvard Business Review 81 (2003), no. 10, 14--16.

\bibitem{Bac}
L. Bacon, Once and for all: tech is not a meritocracy, Quartz, March 27, 2013. http://qz.com/66866/once-and-for-all-tech-is-not-a-meritocracy (accessed July 27, 2014)

\bibitem{Bar}
B. A. Barres, Does gender matter? Nature 442 (2006), 133--136.

\bibitem{BGRL}
S. L. Beilock, E. A. Gunderson, G. Ramirez, and S. C. Levine, Female teachers' math anxiety affects girls' math achievement, Proceedings of the National Academy of Sciences of the USA 107 (2010), no. 5, 1860--1863.

\bibitem{BMN}
M. Biernat, M. Manis, and T. E. Nelson, Stereotypes and standards of judgment, Journal of Personality and Social Psychology 60 (1991), no. 4, 485--499.

\bibitem{BBL}
H. R. Bowles, L. Babcock, and L. Lai, Social incentives for gender differences in the propensity to initiate negotiations: sometimes it does hurt to ask, Organizational Behavior and Human Decision Processes 103 (2007), 84--103.

\bibitem{BGS}
V. L. Brescoll, J. Glass, and A. Sedlovskaya, Ask and ye shall receive? the dynamics of employer-provided flexible work options and the need for public policy, Journal of Social Issues 69 (2013), no. 2, 367--388.

\bibitem{BHKM}
A. W. Brooks, L. Huang, S. W. Kearney, and F. E. Murray, Investors prefer entrepreneurial ventures pitched by attractive men, Proceedings of the National Academy of Science of the USA 111 (2014), no. 12, 4427--4431.

\bibitem{CE}
L. L. Carli and A. H. Eagly, Gender, hierarchy, and leadership: an introduction, Journal of Social Issues 57 (2001), no. 4, 629--636.

\bibitem{CGFSH}
M. Carnes, S. Geller, E. Fine, J. Sheridan, and J. Handelsman, NIH Director's Pioneer Awards: could the selection process be biased against women?, J. Womens Health 14 (2005), no. 8, 684--691.

\bibitem{CH}
A. Casadevall and J. Handelsman, The presence of female conveners correlates with a higher proportion of female speakers at scientific symposia, mBio 5 (2014), no. 1, e00846-13.

\bibitem{CI}
P. R. Clance and S. Imes, The imposter phenomenon in high achieving women: dynamics and therapeutic intervention, Psychotherapy: Theory, Research, and Practice 15 (1978), no. 3, 241--247.

\bibitem{CMR}
R. Cleary, J. W. Maxwell, and C. Rose, Fall 2012 departmental profile report, Notices of the AMS 61 (2014), no. 2, 158--167.

\bibitem{CBP}
S. J. Correll, S. Benard, and I. Paik, Getting a job: is there a motherhood penalty?, American Journal of Sociology 112 (2007), no. 5, 1297--1339.

\bibitem{DK}
D. J. Dean and J. B. Koster, Equitable Solutions for Retaining a Robust STEM Workforce: Beyond best practices, Academic Press, 2014.

\bibitem{DPAHV}
P. G. Devine, E. A. Plant, D. M. Amodio, E. Harmon-Jones, and S. L. Vance, The regulation of explicit and implicit race bias the role of motivations to respond without prejudice, Journal of Personality and Social Psychology 82 (2002), no. 5, 835--848.

\bibitem{DE}
T. A. DiPrete and G. M. Eirich, Cumulative advantage as a mechanism for inequality: a review of theoretical and empirical developments, Annual Review of Sociology 32 (2006), 271--297.

\bibitem{EK}
A. H. Eagly and S. J. Karau, Role congruity theory of prejudice toward female leaders, Psychological Review 109 (2002), no. 3, 573--598.

\bibitem{EHL}
N. M. Else-Quest, J. S. Hyde, and M. C. Linn, Cross-national patterns of gender differences in mathematics: a meta-analysis, Psychological Bulletin 136 (2010), no. 1, 103--127.

\bibitem{EIK}
R. J. Ely, H. Ibarra, and D. M. Kolb, Taking gender into account: theory and design for women's leadership development programs, Academy of Management Learning \& Education 10 (2011), no. 3, 474--493.

\bibitem{FP:C}
Feminist Philosophers, Childcare at conferences: how to do it, August 23, 2010; Childcare at conferences: how to do it (2), September 1, 2010. http://feministphilosophers.wordpress.com/2010/08/23/childcare-at-conferences-how-to-do-it and http://feministphilosophers.wordpress.com/2010/09/01/childcare-at-conferences-how-to-do-it-2 (accessed July 29, 2014)

\bibitem{FP:G}
Feminist Philosophers, Gendered conference campaign, initial version December 9, 2009. http://feministphilosophers.wordpress.com/gendered-conference-campaign (accessed July 28, 2014)

\bibitem{F}
Feministe, Female conference speaker bingo: a bingo card full of excuses for not having more female speakers at STEM conferences. http://www.feministe.us/blog/archives/2012/09/24/why-arent-there-more-women-at-stem-conferences-this-time-its-statistical/female-conference-speaker-bingo (accessed July 28, 2014)

\bibitem{GF:C}
Geek Feminism, Conference anti-harassment/Policy. http://geekfeminism.wikia.com/wiki/Conference\_anti-harassment/Policy (accessed July 29, 2014)

\bibitem{GF:T}
Geek Feminism, Ten tips for getting more female speakers, August 11, 2009. http://geekfeminism.org/2009/08/11/ten-tips-for-getting-more-women-speaker (accessed July 29, 2014)

\bibitem{G}
A. Gheaus, Three cheers for the token woman, Social Science Research Network, March 5, 2013. http://ssrn.com/abstract=2228632 (accessed July 28, 2014)

\bibitem{GR}
C. Goldin and C. Rouse, Orchestrating impartiality: the impact of ``blind'' auditions on female musicians, The American Economic Review 90 (2000), no. 4, 715--741.

\bibitem{GRD}
C. Good, A. Rattan, and C. S. Dweck, Why do women opt out? Sense of belonging and women's representation of mathematics, Journal of Personality and Social Psychology 102 (2012), no. 4, 700--717.

\bibitem{GBN}
T. Greenwald, M. Banaji, and B. Nosek, Take a test: preliminary information, Project Implicit. \\
https://implicit.harvard.edu/implicit/takeatest.html (accessed August 23, 2014)

\bibitem{GMSZ}
L. Guiso, F. Monte, P. Sapienza, and L. Zingales, Culture, gender, and math, Science 320 (2008), 1164--1165.

\bibitem{HWHH}
A. Hegewisch, C. Williams, H. Hartmann, and S. K. Hudiburg, The gender wage gap: 2013, Institute for Women's Policy Research Fact Sheet \#C423, September 2014. http://www.iwpr.org/publications/pubs/the-gender-wage-gap-2013 (accessed September 16, 2014)

\bibitem{H}
M. E. Heilman, The impact of situational factors on personnel decisions concerning women: varying the sex composition of the applicant pool, Organizational Behavior and Human Performance 26 (1980), 386--395.

\bibitem{HSR}
M. E. Heilman, M. C. Simon, and D. P. Repper, Intentionally favored, unintentionally harmed? impact of sex-based preferential selection on self-perceptions and self-evaluations, Journal of Applied Psychology 72 (1987), no. 1, 62--68.

\bibitem{HWFT}
M. E. Heilman, A. S. Wallen, D. Fuchs, and M. M. Tamkins, Penalties for success: reactions to women who succeed at male gender-typed tasks, Journal of Applied Psychology 89 (2004), no. 3, 416--427.

\bibitem{Hess}
A. Hess, Why women aren't welcome on the internet, Pacific Standard 11 (Jan/Feb 2014), 36--47. http://www.psmag.com/navigation/health-and-behavior/women-arent-welcome-internet-72170 (accessed September 8, 2014)

\bibitem{HK}
D. E. Ho and M. G. Kelman, "Does class size affect the gender gap? a natural experiment in law", Journal of Legal Studies 43 (2014), 291--321.

\bibitem{HoMa}
A. Hochschild and A. Machung, The Second Shift: Working families and the revolution at home, Penguin Books, 2012 (revised).

\bibitem{HLLEW}
J. S. Hyde, S. M. Lindberg, M. C. Linn, A. Ellis, and C. Williams, Gender similarities characterize math performance, Science 321 (2008), 494--495.

\bibitem{HM}
J. S. Hyde and J. E. Mertz, Gender, culture, and mathematics performance, Proc. Nat. Acad. Sci. USA 106 (2009), 8801--8807.

\bibitem{IMA}
Institute for Mathematics and its Applications, Diversity at IMA. https://www.ima.umn.edu/diversity (accessed August 23, 2014)

\bibitem{IYH}
L. A. Isbell, T. P. Young, and A. H. Harcourt, Stag parties linger: continued gender bias in a female-rich scientific discipline, PLOS ONE 7 (2012), no. 11, e49682.

\bibitem{JMM}
Joint Mathematics Meetings, Child care grants, 2015 Joint Mathematics Meetings (San Antonio, TX). http://jointmathematicsmeetings.org/meetings/national/jmm2015/2168\_childcare (accessed September 16, 2014)

\bibitem{JSConf}
JSConf EU 2012, Beating the odds---how we got 25\% women speakers for JSConf EU 2012. http://2012.jsconf.eu/2012/09/17/beating-the-odds-how-we-got-25-percent-women-speakers.html (accessed July 29, 2014)

\bibitem{KM}
J. M. Kane and J. E. Mertz, Debunking myths about gender and mathematics performance, Notices of the AMS 59 (2012), no. 1, 10--21.

\bibitem{Kap}
K. Kaplan, Unmasking the impostor, Nature 459 (2009), 468--469.

\bibitem{Kas}
E. Kaschak, Sex bias in student evaluations of college professors, Psychology of Women Quarterly 2 (1978), no. 3, 235--242.

\bibitem{KGH}
S. Knobloch-Westerwick, C. J. Glynn, and M. Huge, The Matilda effect in science communication: an experiment on gender bias in publication quality perceptions and collaboration interest, Science Communication 35 (2013), no. 5, 603--625.

\bibitem{Keyfitz}
B. L. Keyfitz {\em et al.}, Women mathematicians in the academic ranks: a call to action, BIRS Workshop on Women in Mathematics (September 2006), final report. \\
http://www.birs.ca/events/2006/5-day-workshops/06w5504 (accessed November 17, 2014)

\bibitem{Lalin}
M. Lalin, Attending conferences with small children, {\em What's new}, August 20, 2014. http://terrytao.wordpress.com/2014/08/20/matilde-lalin-attending-conferences-with-small-children (accessed August 20, 2014)

\bibitem{LHPL}
S. M. Lindberg, J. S. Hyde, J. L. Petersen, and M. C. Linn, New trends in gender and mathematics performance: A meta-analysis, Psychol. Bull. 136 (2010), 1123--1135.

\vskip4pt\bibitem{Martin}
G. Martin, An annotated bibliography of work related to gender in science. http://arxiv.org/abs/1412.4104

\bibitem{Medina}
H. A. Medina, Doctorate degrees in mathematics earned by blacks, Hispanics/Latinos, and Native Americans: a look at the numbers, Notices of the AMS 51 (2004), no. 7, 772--775.

\bibitem{MAC}
K. L. Milkman, M. Akinola, and D. Chugh, What happens before? a field experiment exploring how pay and representation differentially shape bias on the pathway into organizations, preprint.

\bibitem{MDBGH}
C. A. Moss-Racusin, J. F. Dovidio, V. L. Brescoll, M. J. Graham, and J. Handelsman, Science faculty's subtle gender biases favor male students, Proceedings of the National Academy of Science of the USA 109 (2012), no. 41, 16474--16479.

\bibitem{Munsch}
C. Munsch, Flexible work, flexible penalties: the effect of gender, childcare, and type of request on the flexibility bias, preprint.

\bibitem{MD}
M. C. Murphy and C. S. Dweck, A culture of genius: How an organization's lay theories shape people's cognition, affect, and behavior, Personality and Social Psychology Bulletin 36 (2010), 283--296.

\bibitem{NSF}
National Science Foundation Mathematical Sciences Institutes, Joint math institutes expertise database. \\ http://www.mathinstitutes.org/diversity\_database (accessed September 15, 2014)

\bibitem{N}
I. Neath, How to improve your teaching evaluations without improving your teaching, Psychological Reports 78 (1996), 1363--1372.

\bibitem{Newsome}
J. L. Newsome, The Chemistry PhD: the impact on women's retention, Royal Society of Chemistry/UK Resource Centre for Women in SET. \\ http://www.biochemistry.org/Portals/0/SciencePolicy/Docs/Chemistry\%20Report\%20For\%20Web.pdf

\bibitem{NV}
M. Niederle and L. Vesterlund, Explaining the gender gap in math test scores: the role of competition, Journal of Economic Perspectives 24 (2010), no. 2, 129--144.

\bibitem{OM}
O'Reilly Media, Conference diversity. http://cdn.oreillystatic.com/en/assets/1/eventprovider/1/ConfDiversity.pdf (accessed July 29, 2014)

\bibitem{OECD}
Organisation for Economic Co-operation and Development, How do girls compare to boys in mathematics skills?, PISA 2009 at a Glance, OECD Publishing, 2011.

\bibitem{P}
A. Prasad, Conference diversity distribution calculator. http://aanandprasad.com/diversity-calculator

\bibitem{Polymath}
D. H. J. Polymath, The ``bounded gaps between primes'' polymath project---a retrospective, preprint. arXiv: 1409.8361

\bibitem{PGR}
E. Pronin, T. Gilovich, and L. Ross, Objectivity in the eye of the beholder: perceptions of bias in self versus others, Psychological Review 111 (2004), 781--799.

\bibitem{RSZ}
E. Reuben, P. Sapienza, and L. Zingales, How stereotypes impair women's careers in science, Proceedings of the National Academy of Science of the USA 111 (2014), no. 12, 4403--4408.

\bibitem{R}
C. L. Ridgeway, Gender, status, and leadership, Journal of Social Issues 57 (2001), no. 4, 637--655.

\bibitem{RC}
C. L. Ridgeway and S. J. Correll, Unpacking the gender system: a theoretical perspective on gender beliefs and social relations, Gender and Society 18 (2004), no. 4, 510--531.

\bibitem{Ries}
E. Ries, Why diversity matters (the meritocracy business), Startup Lessons Learned, February 22, 2010. \\
http://www.startuplessonslearned.com/2010/02/why-diversity-matter-meritocracy.html (accessed August 16, 2014)

\bibitem{RM}
E. Ries and S. Milstein, Seeking speakers, Startup Lessons Learned, August 8, 2012. \\
http://www.startuplessonslearned.com/2012/08/seeking-speakers.html (accessed July 29, 2014)

\bibitem{RP}
S. A. Rogier and M. Y. Padgett, The impact of utilizing a flexible work schedule on the perceived career advancement potential of women, Human Resource Development Quarterly 15 (2004), no. 1, 89--106.

\bibitem{Sa}
E. Sander, AWM childcare statement, Association for Women in Mathematics, November 2010. http://sites.google.com/site/awmmath/awm-resources/policy-and-advocacy/awm-childcare-statement (accessed July 29, 2014)

\bibitem{SS}
J. M. Sheltzer and J. C. Smith, Elite male faculty in the life sciences employ fewer women, Proceedings of the National Academy of Science of the USA 111 (2014), no. 28, 10107--10112.

\bibitem{SK}
L. Sinclair and Z. Kunda, Motivated stereotyping of women: she's fine if she praised me but incompetent if she criticized me, Personality and Social Psychology Bulletin 25 (2000), no. 11, 1329--1342.

\bibitem{Snyder}
K. Snyder, The abrasiveness trap: high-achieving men and women are described differently in reviews, Fortune.com, August 26, 2014. http://fortune.com/2014/08/26/performance-review-gender-bias (accessed August 28, 2014)

\bibitem{SSQ}
S. J. Spencer, C. M. Steele, and D. M. Quinn, Stereotype threat and women's math performance, Journal of Experimental Social Psychology 35 (1999), 4--28.

\bibitem{St}
C. Stanton, How I got 50\% women speakers at my tech conference, Geek Feminism, May 21, 2012. http://geekfeminism.org/2012/05/21/how-i-got-50-women-speakers-at-my-tech-conference (accessed July 29, 2014)

\bibitem{SAR}
R. E. Steinpreis, K. A. Anders, and D. Ritzke, The impact of gender on the review of the curricula vitae of job applicants and tenure candidates: a national empirical study, Sex Roles 41 (1999), no. 7/8, 509--528.

\bibitem{Su}
J. Surowiecki, The difference difference makes: waggle dances, the Bay of Pigs, and the value of diversity, in The Wisdom of Crowds, Doubleday, 2004, 23--39.

\bibitem{TP}
F. Trix and C. Psenka, Exploring the color of glass: letters of recommendation for female and male medical faculty, Discourse \& Society 14 (2003), no. 2, 191--220.

\bibitem{Tyson}
N. D. Tyson, response to question during panel discussion, The Secular Society and its Enemies, Center for Inquiry, New York, 2007. response http://www.youtube.com/watch?v=z7ihNLEDiuM ; conference web site http://www.centerforinquiry.net/secularsociety (both accessed August 16, 2014)

\bibitem{UC}
E. L. Uhlmann and G. L. Cohen, ``'I think it, therefore it's true': effects of self-perceived objectivity on hiring discrimination'', Organizational Behavior and Human Decision Processes 104 (2007), 207--223.

\bibitem{U}
D. H. Uttal, Beliefs about genetic influences on mathematics achievement: A cross-cultural comparison, Genetica 99 (1997), 165--172.

\bibitem{Valian}
V. Valian, Recruitment and retention: guidelines for chairs, heads, and deans, The Gender Equity Project, Hunter College, City University of New York, updated February 2011. \\
http://www.hunter.cuny.edu/genderequity/resources/equitymaterials (accessed November 17, 2014)

\bibitem{VS}
V. Valian, D. Sperber, et al., For gender equality at academic conferences. \\ http://forgenderequityatconferences.blogspot.fr ; http://www.gopetition.com/petitions/commitment-to-gender-equity-at-scholarly-conferences.html (accessed July 28, 2014)

\bibitem{VMC}
D. van Dijk, O. Manor, and L. B. Carey, Publication metrics and success on the academic job market, Current Biology 24 (2014), no. 11, R516--R517.

\bibitem{W}
K. Wellhousen, Do's and don'ts for eliminating hidden bias, Childhood Education 73 (1996), no. 1, 36--39.

\bibitem{WW}
C. Wenner\r as and A. Wold, Nepotism and sexism in peer-review, Nature 387 (1997), 341--343.

\bibitem{WIN}
Women in Number Theory, Female Number Theorists. \\
http://womeninnumbertheory.org/index.php?option=com\_content\&view=section\&id=6\&Itemid=13 (accessed November 17, 2014)

\bibitem{WISELI:o}
WISELI, online brochures and booklets, Women in Science \& Engineering Leadership Institute (Madison). Advancing women in science and engineering: advice to the top, http://wiseli.engr.wisc.edu/docs/AdviceTopBrochure.pdf ; Benefits and challenges of diversity in academic settings, http://wiseli.engr.wisc.edu/docs/Benefits\_Challenges.pdf ; Fostering success for women in science and engineering, http://wiseli.engr.wisc.edu/docs/FosteringSuccessBrochure.pdf ; Reviewing applicants: research on bias and assumptions, http://wiseli.engr.wisc.edu/docs/BiasBrochure\_3rdEd.pdf . Accessed July 28, 2014.

\bibitem{WISELI:S}
WISELI, Searching for Excellence \& Diversity: A guide for search committees, Women in Science \& Engineering Leadership Institute (Madison), 2012.

\end{thebibliography}
\end{document}